\RequirePackage{cmap}
\documentclass[12pt, a4paper]{article}
\usepackage[utf8]{inputenc}
\usepackage[T2A]{fontenc}
\usepackage[a4paper,hmargin=2cm,vmargin=2cm]{geometry}
\usepackage[english]{babel}
\usepackage{indentfirst, enumitem}

\usepackage{latexsym} 
\usepackage{amsfonts} 

\usepackage{xcolor}
\usepackage{hyperref}[]
\definecolor{linkcolor}{HTML}{600000}
\definecolor{urlcolor}{HTML}{000050}
\definecolor{citecolor}{HTML}{006000}
\hypersetup{linkcolor=linkcolor,urlcolor=urlcolor, citecolor=citecolor, colorlinks=true}

\ifx\pdfoutput\undefined
\usepackage{graphicx}
\else
\usepackage[pdftex]{graphicx}
\fi
\usepackage{wrapfig}

\usepackage{amssymb}
\usepackage{amsthm}
\usepackage{amsmath}

\DeclareMathOperator{\lk}{lk}
\newtheorem{theorem}{Theorem}[section]

\newtheorem{lemm}[theorem]{Lemma}

\newtheorem{examp}[theorem]{Example}

\newtheorem{conj}[theorem]{Conjecture}
\theoremstyle{definition}
\newtheorem{defin}[theorem]{Definition}
\theoremstyle{remark}
\newtheorem{remark}[theorem]{Remark}

\newcommand{\R}{\mathbb R}

\newcommand*{\hm}[1]{#1\nobreak\discretionary{}%
	{\hbox{$\mathsurround=0pt #1$}}{}}
\newcommand*{\hide}[1]{}
\begin{document}

\author{Mikhail Fedorov\footnote{fedorov.mikhail.s@gmail.com, National Research University Higher School of Economics, Faculty of Mathematics, author is supported by RFBR grant 19-01-00169}}
\title{A description of values of Seifert form for\\punctured n-manifolds in (2n-1)-space.}

\date{}
\maketitle
\begin{abstract}
 We study Seifert linking form which is an invariant of embeddings of punctured $n$-manifolds in $\R^{2n-1}$. 
 For punctured $n$-manifold $N_0$ the values of this invariant are integer valued bilinear symmetric forms on $H_{n-1}(N_0;\mathbb Z)$. We prove that value modulo two of this invariant at $x, y \in H_{n-1}(N_0;\mathbb Z)$ equals $\mathrm{PD}\bar w_{n-2}(N_0)\cap\rho_2x\cap\rho_2y$, where $\mathrm{PD}\bar w_{n-2}(N_0)$ is Poincare dual to Steifel-Whitney class.
 We also prove that any such form can be realized by some embedding $N_0\to\mathbb R^{2n-1}$.

 Also, we survey known results on classification of embeddings of connected manifolds with non-empty boundary.
\end{abstract}
\tableofcontents

\section{Introduction}\label{sec::intro}

Classification of embeddings up to isotopy is one of the basic problems in topology. For a textbook on embeddings see \cite[\S6]{Wa16}, for a survey on classification problem see \cite{Sk16}.
In this paper we explore embeddings of $n$-manifolds with non-empty boundary in $\R^{2n-1}$.\footnote{
Part of this work coincides with my previous term paper and is presented for completeness.
Namely, we repeat statements (but not proofs) of Theorems~\ref{thm::simple_embeddable}, \ref{thm::simple_unknoting}, \ref{thm::k_connect_embeddable}, \ref{thm::k_connect_unknot}, we repeat Definition~\ref{def::linking_form}, Example~\ref{exm::linked_boundary}(a), Lemmas~\ref{lmm::normal_field_exists} and \ref{lmm::L_equality}.
}

This is interesting because results known before unpublished paper \cite{To10} do not
describe this case. For example Theorem~\ref{thm::embeddings_as_immersions} for
embeddings of punctured $k$-connected manifolds into $\R^{2n-k-1}$ for $k>0$, fails for
$k=0$. In \cite{To10, To12} there appears a classification that is based on
\href{http://www.map.mpim-bonn.mpg.de/Some_calculations_involving_configuration_spaces_of_distinct_points}{the
Haefliger-Weber deleted product criterion} \cite[\S5]{Sk06}. My initial aim was to bring
Tonkonog's results to publishable form. However, a more direct approach to
classification was discovered, see Lemma~\ref{lmm::twisting}. We conjecture that there
is a similar proof of injectivity using handle decomposition of a manifold.

We study in \S\S\ref{sec::linking_form}, \ref{sec::descript_values} a non-trivial generalization of Seifert linking form of $2$-manifolds with non-empty boundary in $\mathbb R^3$ to an invariant of punctured $n$-manifold in $\R^{2n-1}$.
This generalization is proposed by D.~Tonkonog in \cite{To12} and is based on \cite{Sa99}, see Definition~\ref{def::linking_form}. Similar invariant is studied in \cite{RWZ}.
The main result of this paper is the description of realizable values of this invariant, see Theorem~\ref{thm::L_is_surjective}.
We give a more direct proof of Theorem~\ref{thm::L_is_surjective} than in \cite{To10}.
Conjecture~\ref{conj::punctured_classification} states that this invariant is injective.
Note that \S\S\ref{sec::linking_form}, \ref{sec::descript_values} are independent of previous sections.\footnote{
This part of the current paper uses D.~Tonkonog's draft \cite{To12}.
We borrowed
Example~\ref{exm::Tonkonog_s} \cite[Example 2.3]{To12},
 Figure~\ref{fig:l3dim} \cite{To12},
Definition~\ref{def::linking_form} \cite[Definition 1.2]{To12},
the first bullet point of Lemma~\ref{lmm::L_well_defined} \cite[Lemma 1.2]{To12},
Lemma~\ref{lmm::saeki} \cite[Lemma 1.9]{To12},
Lemma~\ref{lmm::L_equality} \cite[Lemma 1.5]{To12}
and
Conjecture~\ref{conj::punctured_classification} \cite[Theorem 1.3]{To12}
.
}

The classification of embeddings is known to be hard.
To the best of the author's knowledge, at the time of writing only a few classification results for connected manifolds with non-empty boundary in codimension $3$ are known.
We present a survey of known results (except \cite{To10}) on classification of embeddings of connected manifolds with non-empty boundary, see $\S$\ref{sec::general_theorems} and $\S$\ref{sec::generalisations} of this paper.
Resent surveys on classification of embeddings cover only closed manifolds case, see \cite{Sk06}, \cite{RS99}.
For a more detailed survey on classification of embeddings of manifolds with boundary see \href{https://www.map.mpim-bonn.mpg.de/Embeddings_of_manifolds_with_boundary:_classification}{internet version of \S\S\ref{sec::general_theorems},\ref{sec::example} and \ref{sec::generalisations}} \cite{MAb}.

In \S\ref{sec::example} we give several examples of embeddings and their classification.

Throughout the paper if the category is omitted, then we assume the smooth (DIFF) category. Unless otherwise indicated, the word `isotopy' means `ambient isotopy', see \cite{MAi}.
For any manifold $M$ denote by $\mathrm{Emb}^{m}M$ the set of embeddings of $M$ into $\mathbb R^{m}$ up to isotopy.
Denote by $V_{m, k}$ the Steifel manifold of all $k$-frames in $\mathbb R^m$.

For any closed $n$-manifold $N$ we denote the complement in $N$ to an open $n$-ball by $N_0$.
Thus, $N_0$ is a compact $n$-manifold with non-empty boundary and $\partial N_0$ is the $(n-1)$-sphere.
Note that $H_1(N_0, \partial N_0; *)\hm\cong H_1(N; *)$ because $\partial N_0=S^{n-1}$ is $1$-connected for $n\geq2$. Hereinafter we use everywhere $H_1(N; *)$ instead of $H_1(N_0,\partial N_0; *)$.

Unless otherwise indicated the word ``cycle'' means ``an integer simplicial cycle of some triangulation'', see \cite[\S62]{ST80}.
For any two disjoint $a$-cycle $A$ and $(m-a-1)$-cycle $B$ in $\mathbb R^{m}$ define \textbf{the linking coefficient} by $\lk(A, B) = A \cap \delta B$, where $\delta B$ is a general position $(m-a)$-chain such that $\partial \delta B=B$.
See also \cite{MAif}, \cite[\S3 The linking coefficient]{Sk16a}.

\section{Embedding and unknotting theorems} \label{sec::general_theorems}
\begin{theorem} \label{thm::simple_embeddable}
Assume that $N$ is a compact connected $n$-manifold.
\begin{enumerate}[label=(\alph*), noitemsep, nosep]
\item Then $N$ embeds into $\mathbb R^{2n}$.
\item If $N$ has non-empty boundary, then $N$ embeds into $\mathbb R^{2n-1}$.
\end{enumerate}
\end{theorem}
Part (a) is well-known \href{https://en.wikipedia.org/wiki/Whitney_embedding_theorem}{strong Whitney embedding Theorem}.

\begin{proof}[Proof of part (b)]
By \href{https://en.wikipedia.org/wiki/Whitney_immersion_theorem}{strong Whitney
immersion Theorem} there exists an immersion $g\colon N\hm\to\mathbb R^{2n-1}$. Since
$N$ is connected and has non-empty boundary, it follows that $N$ collapses to an
$(n-1)$-dimensional subcomplex $X\subset N$ of some triangulation of $N$. By general
position we may assume that $g|_{X}$ is an embedding, because $2(n-1) < 2n-1$. Since $g$
is an immersion, it follows that $X$ has a sufficiently small regular neighborhood
$M\supset X$ such that $g|_{M}$ is embedding. Since regular neighborhood is unique up to
homeomorphism, there exists a homeomorphism $h\colon N\to M$.
The composition $g\circ h$ is an embedding of $N$.
\end{proof}
This proof is essentially contained in \cite[Theorem 4.6]{Hi61}.
\begin{theorem}\label{thm::simple_unknoting}
Assume that $N$ is a compact connected $n$-manifold and either
\begin{enumerate}[label=(\alph*), noitemsep, nosep]
\item $m \ge 2n+1 \ge 5$ or
\item $N$ has non-empty boundary and $m \ge 2n$.
\end{enumerate}
Then any two embeddings of $N$ into $\mathbb R^m$ are isotopic.
\end{theorem}

All theorems stated in \S\S \ref{sec::general_theorems}, \ref{sec::generalisations} for
embeddings of manifolds with non-empty boundary can be proved using analogous results for
immersions of manifolds, general position and
\href{http://www.map.mpim-bonn.mpg.de/Isotopy}{Concordance Implies Isotopy Theorem}
\cite[\S2]{MAi}. We give proofs for Theorem~\ref{thm::simple_embeddable}(b) and
Theorem~\ref{thm::k_connect_unknot}(b) to illustrate these ideas.

Theorem~\ref{thm::simple_unknoting}(a) is the Whitney-Wu Unknotting Theorem, see \cite[Theorems 2.1, 2.2]{Sk16}.

Theorem~\ref{thm::simple_unknoting}(b) for $n>2$ can be found%
\footnote{
We do not claim the references we give here and in \S\ref{sec::generalisations} are references to original proofs.
}
in \cite[$\S$ 4, Corollary 5]{Ed68}. Case $n=1$ is clear. Case $n=2$ can be proved using the above ideas.

\section{Examples of non-isotopic embeddings}\label{sec::example}
Examples we give in this section are simple, so we assume that they are folklore.

Denote $1_k:=(0,\ldots,0,1)\in S^k$.
\begin{examp} \label{exm::linked_boundary}
Let $N=S^k\times [0, 1]$ be the cylinder over $S^k$.
\begin{enumerate}[label=(\alph*), noitemsep, nosep]
\item There exist non-isotopic embeddings of $N$ into $\mathbb R^{2k+1}$.

\item Moreover, for each $a\in\mathbb Z$ there exist an embedding  $f\colon N\to\mathbb R^{2k+1}$ such that $\lambda(f) \hm{:=}\mathrm{lk}(f(S^k\times 0), f(S^k\times 1))=a$.

\item Moreover, $\lambda\colon \mathrm{Emb}^{2k+1}N\to\mathbb Z$ is well-defined and is a 1--1 correspondence for $k\geqslant2$.
\end{enumerate}
\end{examp}

\begin{proof}[Proof of part (b)]  For $k=1$ the cylinder is a ribbon.
Informally speaking by twisting a ribbon one can obtain an arbitrary value of the linking coefficient.

Now we give a rigorous construction.
Let $h\colon S^k\to S^k$ be a map of degree $a$.
(To prove part (a) it is sufficient to repeat this construction twice
with the identity map of $S^k$ as a map of degree one
and with the constant map as a map of degree zero.)
Define an embedding $g\colon S^k\times [0, 1] \to D^{k+1}\times S^k$ by the formula $g(x, t) = (h(x)t, x)$.

Let $f=\mathrm i\circ g$, where $\mathrm i = \mathrm i_{2k+1, k}\colon D^{k+1}\times S^k \to \mathbb R^{2k+1}$ is \href{http://www.map.mpim-bonn.mpg.de/Embeddings_in_Euclidean_space:_an_introduction_to_their_classification#Notation_and_conventions}{the standard embedding} \cite[\S 3 Notation and conventions]{Sk16}.

Let us show that $\mathrm{lk}(f(S^k\times0), f(S^k\times1)) = a$.
Consider the standard $k$-disk $D'$ in $\R^{2k+1}$ in the hyperplane $x_1=\cdots=x_{k-1}=0$ such that
$$\partial D' = f(S^k\times{0})
\quad\text{and}\quad
D'\cap \mathrm i(\partial D^{k+1}\times S^k) = \mathrm i(1_k\times S^k).$$
The linking coefficient of $f(S^k\times0)$ and $f(S^k\times1)$ is the sum of signs of intersection points in $f(S^k\times{0})\cap D'=\{(x, 1)\mid h(x) = 1_k\}$.
This sum equals $\deg h = a$.
\end{proof}
\begin{proof}[Proof of part (c)] Clearly $\lambda$ is well-defined.
By (b) $\lambda$ is surjective.
Now take any two embeddings $f_1, f_2$ such that $\lambda([f_1]) = \lambda([f_2])$.
Each embedding of a cylinder gives an embedding of a sphere with a (nowhere vanishing) normal field.
Since $k\geqslant 2$, \href{http://www.map.mpim-bonn.mpg.de/Embeddings_in_Euclidean_space:_an_introduction_to_their_classification#Unknotting_theorems}{the Unknotting Spheres Theorem} \cite[\S2, Teorem 2.3]{Sk16} implies that there exists an isotopy between $f_1|_{S^k\times 0}$ and $f_2|_{S^k\times 0}$.
Thus, we can assume that $f_1|_{S^k\times 0} = f_2|_{S^k\times 0}$.
Since $\lambda([f_1]) = \lambda([f_2])$, it follows that  normal fields on $f_1(S^k\times 0)$ and $f_2(S^k\times 0)$ are homotopic in class of normal fields.
So $f_1$ and $f_2$ are isotopic.
\end{proof}

Recall that we denote the complement in $N$ to an open $n$-ball by $N_0$.

\begin{examp}\label{exm::punctured_torus}
Let $N=S^k\times S^1$. Assume that $k>2$ and that $S^k\times\{1_1, -1_1\}\subset N_0$.
Then there exists a 1--1 correspondence $\lambda\colon\mathrm{Emb}^{2k+1}N_0\to\mathbb Z$ defined by the formula $\lambda([f])=\mathrm{lk}(f(S^k\hm\times1_1), f(S^k\hm\times-1_1))$.
\end{examp}
The surjectivity of $\lambda$ is proved analogously to Example \ref{exm::linked_boundary}(b).
The injectivity of $\lambda$ follows from Example~\ref{exm::linked_boundary}(c)
because the restriction-induced map $\mathrm{Emb}^{2k+1}N_0\to\mathrm{Emb}^{2k+1}(S^k\hm\times[0,1])$ is a 1--1 correspondence by general position.

\begin{examp}\label{exm::sum_of_tori}
Let $N=S^k_a\times S^1 \# S^k_b\times S^1$ be the connected sum of two tori.
\begin{enumerate}[label=(\alph*), noitemsep, nosep]
 \item Then there exists a surjection $\lambda\colon\mathrm{Emb}^{2k+1}N_0\to\mathbb Z$ defined by the formula $\lambda([f])\hm=\mathrm{lk}(f(S^k_a\times1_1), f(S^k_b\times1_1))$.
 \item Then $\lambda\colon \mathrm{Emb}^{2k+1}N_0\to \mathbb Z$ is not injective.
 \item The map $\mathrm{Emb}^{2k+1}N_0\to \mathbb Z^3$ given by
 $$\lk(f(S^k_a\times1_1), f(S^k_b\times1_1)),
 \quad
 \lk(f(S^k_a\times1_1), f(S^k_a\times-1_1))
 \quad\text{and}\quad
 \lk(f(S^k_b\times1_1), f(S^k_b\times-1_1)) $$
is surjective.
\end{enumerate}
\end{examp}
To prove part (a) it is sufficient to take linked $k$-spheres in $\mathbb R^{2k+1}$ with normal fields.
The normal fields induce ribbons containing these two spheres.
Analogously to Example~\ref{exm::punctured_torus} embeddings of ribbons gives embeddings of punctured tori.
Consider an embedded boundary connected sum of these two tori.
Part (b) follows from Example~\ref{exm::punctured_torus}.
Part (c) is proved analogously to part (a) and surjectivity of $\lambda$ in Example~\ref{exm::punctured_torus}.

In all examples above we chose disjoint submanifolds in $N$ representing two certain homology classes, for example $S^k\times0$ and $S^k\times1$ in Example~\ref{exm::linked_boundary}.
For arbitrary $n$-manifold and homology classes this may be impossible.
For example, in $(S^1)^n$ there exist $(n-1)$-cycles such that any $(n-1)$-submanifolds representing them intersect.
For intersecting submanifolds the linking coefficient is not defined.
However, it is possible to define an analogue of invariant $\lambda$ based on linking coefficient, see the following example and Definition~\ref{def::linking_form}.

\begin{figure}[h]
\centering
\includegraphics{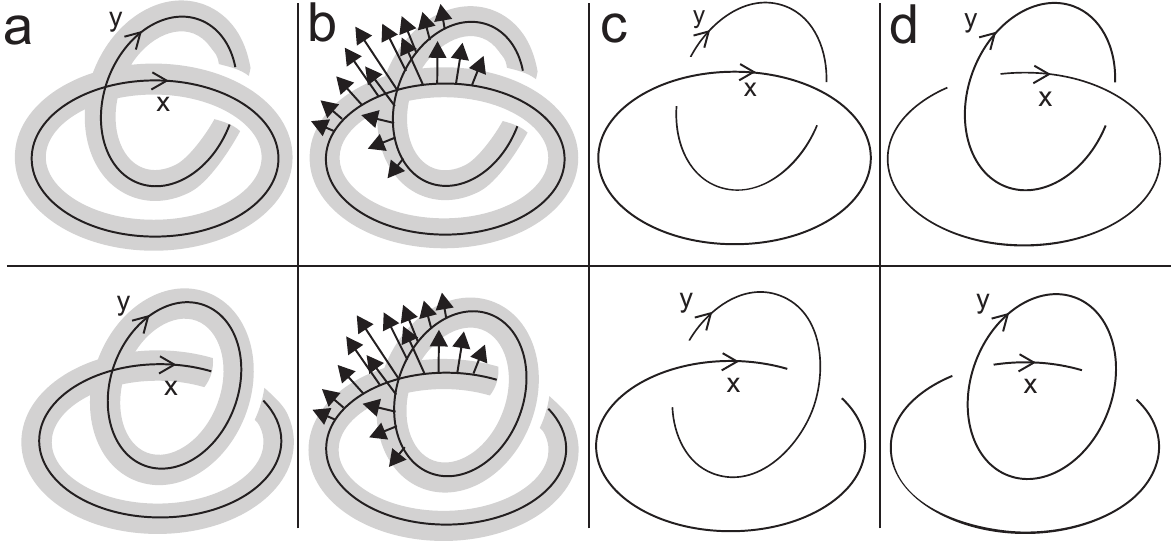}
\caption{
(a) Embeddings $f_1$ (top) and $f_2$ (bottom);\\
(b) the vectors depict the difference $\xi_{f_i}-f_i$, $i=1,2$, so their ends define $\xi_{f_i}$;\\
(c) linking $\xi_{f_i}(x)\sqcup f_i(y)$;\\
(d) linking $\bar \xi_{f_i}(x)\sqcup f_i(y)$, where $\bar \xi_i$ is the normal field opposite to $\xi_i$.
}
\label{fig:l3dim}
\end{figure}

\begin{examp}\label{exm::Tonkonog_s}
(a) Let $N_0$ be the oriented punctured 2-torus containing the meridian $x$ and the parallel $y$ of the torus.
For each embedding $f\colon N_0\to\mathbb R^3$ denote by $\xi_f\colon N_0\hm\to\mathbb R^3$ the normal field on $f(N_0)$ defined by the orientations on $N_0$ and on $\mathbb R^3$ (see figure \ref{fig:l3dim}(b)).
Then the map $\lambda\colon\mathrm{Emb}^3 N_0\to\mathbb Z$ defined by the formula $\lambda([f])=\mathrm{lk}(\xi_f(x), f(y))$ is a surjection.

(b) Let $f_1,f_2\colon N_0\to\R^3$ be two embeddings shown in figure \ref{fig:l3dim}(a).
Figure \ref{fig:l3dim}(c) shows that $\lambda(f_1)=1$ and $\lambda(f_2)=0$ which proves that $f_1$ and $f_2$ are not isotopic.
(Notice that the restrictions of $f_1$ and $f_2$ to $x\cup y$ are isotopic!).
If in definition of $\lambda$ we change $\xi_f$ to the opposite normal vector field $\bar\xi_f$ then values $\lambda(f_1)$ and $\lambda(f_2)$ will change but will still be different (see figure \ref{fig:l3dim}(d)).
\end{examp}

\section{A generalization to highly-connected manifolds} \label{sec::generalisations}

In this section we do not consider manifolds with arbitrary non-empty boundary. Instead, we consider punctured manifolds to simplify statements. However, analogous results seem to hold for manifolds with arbitrary boundary, see \cite[\S7 Comments on non-spherical boundary]{MAb}.

\begin{theorem}\label{thm::k_connect_embeddable}
Assume that $N$ is a closed $k$-connected $n$-manifold.
\begin{enumerate}[label=(\alph*), noitemsep, nosep]
\item If $n\geq2k+3$, then $N$ embeds into $\mathbb R^{2n-k}$.
\item $N_0$ embeds into $\mathbb R^{2n-k-1}$.
\end{enumerate}
\end{theorem}
Part (a) of this result is \cite[Existence Theorem (a)]{Ha61}.
The PL case of this result is well-known Penrose-Whitehead-Zeeman Theorem, see \cite[Theorem 1.1]{PWZ61}.

Part (b) is \cite[Corollary 4.2]{Hi61}.

\begin{theorem}\label{thm::k_connect_unknot}
Assume that $N$ is a closed $k$-connected $n$-manifold.
\begin{enumerate}[label=(\alph*), noitemsep, nosep]
\item If $n \geq 2k + 2$ and $m \geq 2n - k + 1$, then any two embeddings of $N$ into $\mathbb R^m$ are isotopic.
\item If $n\ge 2k+2$, $m\geq2n-k$ and $(m,n,k)\neq(7,4,1)$, then any two embeddings of $N_0$ into $\mathbb R^m$ are isotopic.
\end{enumerate}
\end{theorem}

For Theorem~\ref{thm::k_connect_unknot}(a) see \cite[Theorem 2.4]{Sk16}, or \cite[Corollary 2 of Theorem 24 in Chapter 8]{Ze63} and \cite[Existence Theorem (b) in p. 47]{Ha61}.

For $k>1$ theorem~\ref{thm::k_connect_unknot}(b) is a particular case of Theorem~\ref{thm::embeddings_as_immersions} (but of course can be proved in a direct way).
For $k=0$ theorem~\ref{thm::k_connect_unknot}(b) coincides with Theorem~\ref{thm::simple_unknoting}(b).
For $k=1$ the result seems to be known, but we did not find a reference except unpublished paper \cite{To10} (where this Theorem is proven without exception $(m,n,k)=(7,4,1)$).
We present a proof for which we need the following
\begin{theorem}\label{thm::spine}
Assume $N$ is a closed $k$-connected $n$-manifold and $n\geq k+3$, $(n,k)\hm\notin \{(5,2), (4,1)\}$. Then there exists a $(n-k-1)$-dimensional subcomplex $X\subset N_0$ of some triangulation of $N_0$ such that $N_0$ collapses to $X$.
\end{theorem}
This result can be found in \cite[Theorem 5.5]{Wa64}, see also \cite[Lemma 5.1 and Remark 5.2]{Ho69}. A counterexample for $(n, k) = (4, 1)$ is presented in \cite{LL18}, see also \cite{Sk19}.

Denote by $TN_0$ the tangent bundle of $N_0$.

\begin{proof}[Proof of Theorem~\ref{thm::k_connect_unknot} (b) for $k=1$]
If $(n,k) = (4, 1)$ then $m\geq8=2n$.
So this case is a particular case of Theorem~\ref{thm::simple_unknoting}(b).
Assume $(n,k)\neq(4,1)$.
Smale-Hirsch Theorem \cite[Theorem 6.1]{Hi59} implies that regular homotopy classes of immersions of $N_0$ in $\mathbb R^m$ are in 1--1 correspondence with homotopy classes of linear monomorphisms $TN_0\hm\to\mathbb R^m$.
The obstructions for two linear monomorphisms to be regularly homotopic lie in $H^r(N_0; \pi_r(V_{m, n}))$. All these groups are zero because $N$ is $1$-connected and $\pi_r(V_{m, n})=0$ for $r<n-1$. Hence, any two immersions $N_0\to\mathbb R^m$ are regularly homotopic.

Hence, for every two embeddings $f,g\colon N_0\to\mathbb R^m$ there exists an immersion $F\colon N_0\hm\times[0,1]\hm\to\mathbb R^m\hm\times[0,1]$ such that $F(x, 0) = (f(x), 0)$ and $F(x, 1)=(g(x), 1)$ for each $x\in N_0$.
Theorem~\ref{thm::spine} implies that $N_0$ collapses to an $(n-2)$-dimensional subcomplex $X\subset N_0$ of some triangulation of $N_0$.
By general position we may assume that $F|_{X\times[0,1]}$ is an embedding, because $2(n-1) < m+1$.
Since $F$ is an immersion, it follows that $X$ has a sufficiently small regular neighborhood $M\supset X$ such that $F|_{M\times[0,1]}$ is an embedding.
Since regular neighborhood is unique up to homeomorphism, there exists a homeomorphism $h\colon N_0\to M$.
Thus, the restriction $F|_{M\times[0,1]}$ is a concordance between $f\circ h$ and $g\circ h$.
By \href{http://www.map.mpim-bonn.mpg.de/Isotopy#Concordance}{the Concordance Implies Isotopy Theorem} \cite[\S 2]{MAi} $f\circ h$ and $g\circ h$ are isotopic.
The collapse of $N_0$ to $X$ induce isotopies between $f$ and $f\circ h$, as well as between $g$ and $g\circ h$.
\end{proof}

\begin{theorem} \label{thm::embeddings_as_immersions}
Assume that $N$ is a closed orientable $k$-connected manifold and $(n, k)\hm\notin\{(5,2), (4,1)\}$. 
Then for each $k\ge1$ there exists a 1--1 correspondence
$$
W_0'\colon \mathrm{Emb}^{2n-k-1}N_0\to H_{k+1}(N;\mathbb Z_{(n-k-1)}),
$$
where $W_0'$ denotes the Whitney invariant \cite[Definition 1.3]{Sk10} and $\mathbb Z_{(t)}$ denotes $\mathbb Z$ for $t$ even and $\mathbb Z_2$ for $t$ odd.
\end{theorem}
This theorem was essentially known since 1970s, but presumably was not explicitly stated before \cite[Lemma 2.2($W_0'$)]{Sk10}. The proof below follows \cite[proof of Lemma 2.2($W_0'$)]{Sk10}.
\begin{proof}
Let $Y$ be the set of regular homotopy classes of immersions $N_0\to\mathbb R^{2n-k-1}$. Theorem~\ref{thm::spine} implies that $N_0$ collapses to an $(n-k-1)$-dimensional subcomplex $X\subset N_0$ of some triangulation of $N_0$. So by general position the forgetful map $\mathrm{Emb}^{2n-k-1}N_0\to Y$ is surjective (because $2(n-k-1) < 2n-k-1$) and, for $k\geq1$, injective  (because $2(n-k)<2n-k$). The map $W'_0$ is a composition of the forgetful map and a 1--1 correspondence $\phi\colon Y\to H_{k+1}(N;\mathbb Z_{(n-k-1)})$, which we construct in the next paragraph.

Smale-Hirsh theorem \cite[Theorem 6.1]{Hi59} implies that $Y$ is in 1--1 correspondence with homotopy classes of linear monomorphisms $TN_0\hm\to\mathbb R^{2n-k-1}$. The obstructions for two linear monomorphisms to be homotopic lie in $H^t(N; G_t)$, where $G_t = \pi_t(V_{2n-k-1, n})$. All these groups except group $H^{n-k-1}(N_0;G_{n-k-1})$ are zero because $N$ is $k$-connected and $G_t = 0$ for $t<n-k-1$. We know $G_{n-k-1}\cong \mathbb Z_{(n-k-1)}$. The obstruction is unique, so it gives the 1--1 correspondence $\phi$.
\end{proof}

\begin{remark} \label{rem::fails}
The analogue of Theorem~\ref{thm::embeddings_as_immersions} for $k=0$ fails,
i.e. there is no 1--1 correspondence between $\mathrm{Emb}^{2n-1}N_0$ and $H_1(N;\mathbb Z_2)$.
This is shown by Example~\ref{exm::punctured_torus}.
Moreover, the integer analogue of Theorem~\ref{thm::embeddings_as_immersions} for $k=0$ fails,
i.e. there is no 1--1 correspondence between $\mathrm{Emb}^{2n-1}N_0$ and $H_1(N;\mathbb Z)$.
This is shown by Example~\ref{exm::sum_of_tori}.
For description of $\mathrm{Emb}^{2n-1}N_0$ see Conjecture~\ref{conj::punctured_classification}.
\end{remark}

\section{Seifert linking form} \label{sec::linking_form}

Results stated in \S\S~\ref{sec::general_theorems},~\ref{sec::generalisations} give an impression that embeddings of $N_0$ into $\mathbb R^{m}$ are in 1--1 correspondence to embeddings of $N$ into $\mathbb R^{m+1}$ for any integer $m$. In his unpublished work \cite{To10} D.~Tonkonog shows the first case when $\mathrm{Emb}^mN_0$ and $\mathrm{Emb}^{m+1}N$ are not equal.
We provide an alternative proof of this fact.

Unless otherwise stated, we omit $\mathbb Z$-coefficients from the notation of (co)homology groups.
Throughout this section we assume that $N$ is a closed orientable connected $n$-manifold, $n\geq 4$, $n$ is even,  $H_1(N)$ is torsion free and $f\colon N_0\to\mathbb R^{2n-1}$ is an embedding.

The following folklore result holds.
\begin{lemm}\label{lmm::normal_field_exists}{\normalfont(\cite{To12}, cf. \cite[Lemma 4.1]{Sa99})}
There exists a nowhere vanishing normal vector field on $f(N_0)$.
\end{lemm}
\begin{proof}
There is an obstruction (Euler class) $\bar e=\bar e(f)\in H^{n-1}(N_0)\cong H_1(N)$ to existence of a nowhere vanishing normal vector field on $f(N_0)$.

A normal space to $f(N_0)$ at any point of $f(N_0)$ has dimension $n-1$.
As $n$ is even thus $n-1$ is odd. Thus, if we replace a general position normal field by its opposite then the obstruction will change sign.
Therefore, $\bar e=-\bar e$.
Since $H_1(N)$ is torsion free, it follows that $\bar e=0$.

Since $N_0$ has non-empty boundary, it follows that $N_0$ collapses to an $(n-1)$-dimensional subcomplex $X\subset N_0$ of some triangulation of $N_0$.
The dimension of $X$ equals the dimension of normal space to $f(N_0)$ at any point of $f(N_0)$.
Since $\bar e=0$, it follows that there exists a nowhere vanishing normal vector field on $f(N_0)$.
\end{proof}

\begin{defin}\label{def::linking_form}
For every $(n-1)$-cycles $X,Y$ in $N_0$ define \textbf{Seifert linking form} by the formula
$$L(f)(X,Y) = \mathrm{lk}(f(X), \xi(Y)) + \mathrm{lk}(\xi(X), f(Y)),$$
where $\xi\colon N_0\to\mathbb R^{2n-1}$ is a nowhere vanishing normal field on $f(N_0)$.
Note that $L(f)$ is bilinear because $\lk$ is bilinear. And $L(f)$ is symmetric because $\lk$ is symmetric for $n$ even.
\end{defin}

In \cite{Sa99} Seifert linking form is defined for manifolds with normal field for $n=3$ as the first summand. Saeki shows that the first summand is not an invariant. Presumably it is possible to define an analogous invariant for manifolds with normal field for odd $n$.

\begin{lemm}\label{lmm::L_well_defined}
The integer $L(f)(X, Y)$:
\begin{itemize}[noitemsep, nosep]
    \item is well-defined, i.e. is independent of the choice of  $\xi$,
    \item does not change when $X$ or $Y$ are changed to homologous cycles, and
    \item does not change when $f$ is changed to an isotopic embedding.
\end{itemize}
\end{lemm}

We will need the following sublemma.
\begin{lemm} \label{lmm::saeki}
Let $\xi,\xi'\colon N_0\to\mathbb R^{2n-1}$ be two nowhere vanishing normal vector fields on $f(N_0)$.
Then for any $(n-1)$-cycles $X,Y$ in $N_0$ holds
$$\mathrm{lk}(f(X),\xi(Y))-\mathrm{lk}(f(X),\xi'(Y))=d(\xi,\xi')\cap [X]\cap [Y]$$
where $d(\xi,\xi')\in H_2(N_0)$
is (Poincare dual to) the first obstruction to $\xi,\xi'$ being homotopic in the class of the nowhere vanishing vector fields.
\end{lemm}
Lemma~\ref{lmm::saeki} is proved in \cite[Lemma 2.2]{Sa99} for $n=3$, but the proof is valid in all dimensions \cite{To12}.
\begin{proof}[Proof of Lemma \ref{lmm::L_well_defined}.]
The first bullet point follows because for any nowhere vanishing normal field $\xi'$ on $N_0$ we have:
\begin{equation*}
\begin{aligned}
\mathrm{lk}(f(X),\xi(Y))&+\mathrm{lk}(\xi(X),f(Y))
&-\mathrm{lk}(f(X),\xi'(Y))&-\mathrm{lk}(\xi'(X),f(Y))= \\
\mathrm{lk}(f(X),\xi(Y))&+(-1)^n\,\mathrm{lk}(f(Y),\xi(X))
&-\mathrm{lk}(f(X),\xi'(Y))&-(-1)^n\,\mathrm{lk}(f(Y),\xi'(X))=\\
&&=d(\xi,\xi')\cap [X]\cap [Y]+&(-1)^n \,d(\xi,\xi')\cap [Y]\cap [X]=\\
&&=d(\xi,\xi')\cap [X]\cap [Y]&(1+(-1)^n(-1)^{n-1})=0.
\end{aligned}
\end{equation*}
Here the second equality follows from Lemma~\ref{lmm::saeki}.

Now let us prove the second bullet point. Denote by $X, X'$ two homologous $(n-1)$-cycles in $N_0$.
The image of the homology between $X$ and $X'$ is a $n$-chain $\delta X$ in $f(N_0)$ such that $\partial \delta X = f(X) - f(X')$.
Since $\xi$ is a nowhere vanishing normal field on $f(N_0)$, this implies that the supports of $\xi(Y)$ and $\delta X$ are disjoint. Hence, $\mathrm{lk}(f(X), \xi(Y)) = \mathrm{lk}(f(X'), \xi(Y))$.

Now let us prove the third bullet point.
Since an isotopy of $f$ is a map from $\mathbb R^{2n-1}\times [0, 1]$ to $\mathbb R^{2n-1}\times [0, 1]$, it follows that this isotopy gives an isotopy of the link $f(X)\sqcup \xi(Y)$. Now the third bullet point follows because the linking coefficient is preserved under isotopy.
\end{proof}

Lemma~\ref{lmm::L_well_defined} implies that $L(f)$ generates a bilinear form $H_{n-1}(N_0)\times H_{n-1}(N_0)\to\mathbb Z$ denoted by the same notation.

\section{Values of Seifert linking form}\label{sec::descript_values}
In this section we assume that $N$ is a closed orientable connected $n$-manifold, $n\geq 4$, $n$ is even and $H_1(N)$ is torsion free.
Denote by $\rho_2 \colon H_*(N)\to H_*(N;\mathbb Z_2)$ the reduction modulo $2$.

Recall that the dual to \href{http://www.map.mpim-bonn.mpg.de/Stiefel-Whitney_characteristic_classes}{Stiefel-Whitney class} $\mathrm{PD}\bar w_{n-2}(N_0)\in H_2(N_0; \mathbb Z_2)$ is the class of the cycle on which two general position normal fields on $f(N_0)$ are linearly dependent.
\begin{lemm}[\cite{To12}] \label{lmm::L_equality}
Let $f:N_0\to \mathbb R^{2n-1}$ be an embedding.
Then for every $x, y \hm\in H_{n-1}(N_0)$ the following equality holds:
$$\rho_2L(f)(x, y) = \mathrm{PD}\bar w_{n-2}(N_0)\cap\rho_2x\cap\rho_2y.$$
\end{lemm}
In particular this Lemma shows that the Seifert linking form modulo $2$ is independent of the embedding.

Here we provide the proof with more details than in \cite{To12}.

See also an analogous lemma for closed manifolds in \cite[Lemma 2.2]{CS16}.

\begin{proof}[Proof of Lemma \ref{lmm::L_equality}.]
Let $\xi$ be a normal field on $f(N_0)$. Denote by $\bar\xi$ the normal field on $f(N_0)$ opposite to $\xi$. We have
$$
\begin{aligned}
L(f)(x, y) &\underset{2}{\stackrel{(1)}{\equiv}} \mathrm{lk}(f(x), \xi(y)) - \mathrm{lk}(\xi(x), f(y))
	\stackrel{(2)}= \\&=
\mathrm{lk}(f(x), \xi(y)) - \mathrm{lk}(f(x), \bar\xi(y))
    \stackrel{(3)}=
d(\xi, \bar\xi)\cap x\cap y
	.
\end{aligned}
$$
The congruence (1) is clear.

The equality (2) holds because if we shift the link $\xi(x)\sqcup f(y)$ by $\bar\xi$, we get the link $f(x)\sqcup \bar\xi(y)$ and the linking coefficient will not change after this shift.

The equality (3) follows from Lemma \ref{lmm::saeki}.

Now the lemma follows because
$$
\rho_2 d(\xi, \bar\xi) \stackrel{(*)}= \rho_2 d(\xi', \bar\xi) \stackrel{(**)}= \mathrm{PD}\bar w_{n-2}(N_0),
$$
where $\xi'$ is a general perturbation of $\xi$.

The equality (*) holds because $\xi'$ and $\xi$ are homotopic in the class of nowhere vanishing normal vector fields.

In this paragraph we prove the equality (**). Vector fields $\xi$ and $\xi'$ are in general position, to they are linearly dependent on a 2-cycle modulo 2 which we denote by $c$. We have $[c]=\mathrm{PD}\bar w_{n-2}(N_0)$.
By construction vector fields $\xi'$ and $\bar\xi$ are linearly dependent at some point if and only if they are opposite at this point. Hence, the linear homotopy between $\xi'$ and $\bar\xi$ degenerates on $c$.
Since the linear homotopy between $\xi'$ and $\bar\xi$ is a general position homotopy, it follows that $\rho_2 d(\xi', \bar\xi) = [c]$.
\end{proof}

\begin{theorem}[\cite{To12}]\label{thm::L_is_surjective}
For any symmetric bilinear form $\phi\colon H_{n-1}(N_0)\hm\times H_{n-1}(N_0)\to \mathbb Z$ such that $\rho_2\phi(x,y)= \mathrm{PD}\bar w_{n-2}(N_0)\cap\rho_2x\cap \rho_2y$ there exist an embedding $f\colon N_0\to\R^{2n-1}$ such that $L(f) = \phi$.
\end{theorem}

\begin{conj}[\cite{To12}]\label{conj::punctured_classification}
Let $N$ be a closed connected orientable $n$-manifold with $H_1(N)$ torsion-free, $n\ge 4$, $n$ even. Then Seifert linking form $L$ is an injection.

I.e. for any two embeddings $f, g\colon N_0\to\mathbb R^{2n-1}$ if $L(f) = L(g)$ then $f$ and $g$ are isotopic.
\end{conj}

Theorem~\ref{thm::L_is_surjective} and Conjecture~\ref{conj::punctured_classification} correspond to \cite[Theorem~1.3]{To12}, cf. \cite[Conjecture~1]{To10}.
The proof in \cite{To12} is based on \cite{To10} and significantly uses \href{http://www.map.mpim-bonn.mpg.de/Some_calculations_involving_configuration_spaces_of_distinct_points}{the Haefliger-Weber deleted product criterion}.
Theorem~\ref{thm::L_is_surjective} follows from Lemma~\ref{lmm::twisting} and a well-known Lemma~\ref{lmm::linear_algebra}.

\begin{lemm}\label{lmm::twisting}
For each embedding $f\colon N_0\to \mathbb R^{2n-1}$, class $\tilde s\in H_1(N_0, \partial N_0)$ and integer $l\in\mathbb Z$
there exist an embedding $g\colon N_0\to \mathbb R^{2n-1}$ such that
$$L(g)(x, y) - L(f)(x, y) = 2l(\tilde s\cap x) \cdot(\tilde s\cap y)$$
for any $x, y \in H_{n-1}(N_0)$.
\end{lemm}

The idea is to prove this lemma using ``twisting a ribbon containing $\tilde s$'', as in Example~\ref{exm::linked_boundary}.
The proof below does not mention parametric connected sum \cite[\S2.4]{CS16} but in fact uses that construction.

\begin{lemm}\label{lmm::linear_algebra}
Assume that $V$ is a free $\mathbb Z$-module and $m\colon V\times V\to\mathbb Z$ is a symmetric bilinear form.
Then there exists linear functions $s_1,\ldots,s_n\hm\in V^*$ and coefficient $a_1, \ldots ,a_n\in\mathbb Z$ such than $m(x, y)=\sum_{i=1}^{n}a_is_i(x)s_i(y)$ for each $x, y\in V$.
\end{lemm}

\begin{proof}[Proof of Lemma~\ref{lmm::twisting}]
Take an embedding $s\colon D^1\times D^{n-1}\to N_0$ such that $s(\partial D^1 \times D^{n-1})\hm\subset \partial N_0$ and the relative homology class of $s(D^1\times0)$ is $\tilde s$.
Denote by $\xi_f\colon N_0 \to \mathbb R^{2n-1}$ a normal field on $f(N_0)$.
Denote $1_k=(0,\ldots,0,1)\in S^k$.

Any bundle over disk is trivial, so there exists an $n$-framing on the disk $fs(0\hm\times D^{n-1})$
such that the first vector at $fs(0, a)$ is $d(fs)|_{(0,a)}(1_0)$ and the second vector is $\xi_fs(0,a) \hm- fs(0,a)$.
Using this framing one can construct an embedding $\psi \colon D^n\times D^{n - 1} \hm\to \mathbb R^{2n-1}$ such that
$$
\psi(D^1\times D^{n-1})=fs(D^1\times D^{n-1}),\;
\psi(0\times D^{n-1}) = fs(0\times D^{n-1})\;\text{ and }\;
\psi(1_1\times D^{n-1}) = \xi_fs(0\times D^{n-1}),
$$
(note that $1_1\notin D^1$ but $1_1\in S^1$).

\begin{figure}[!h]
\label{fig::twist}
 \centering
 \includegraphics[scale=0.7]{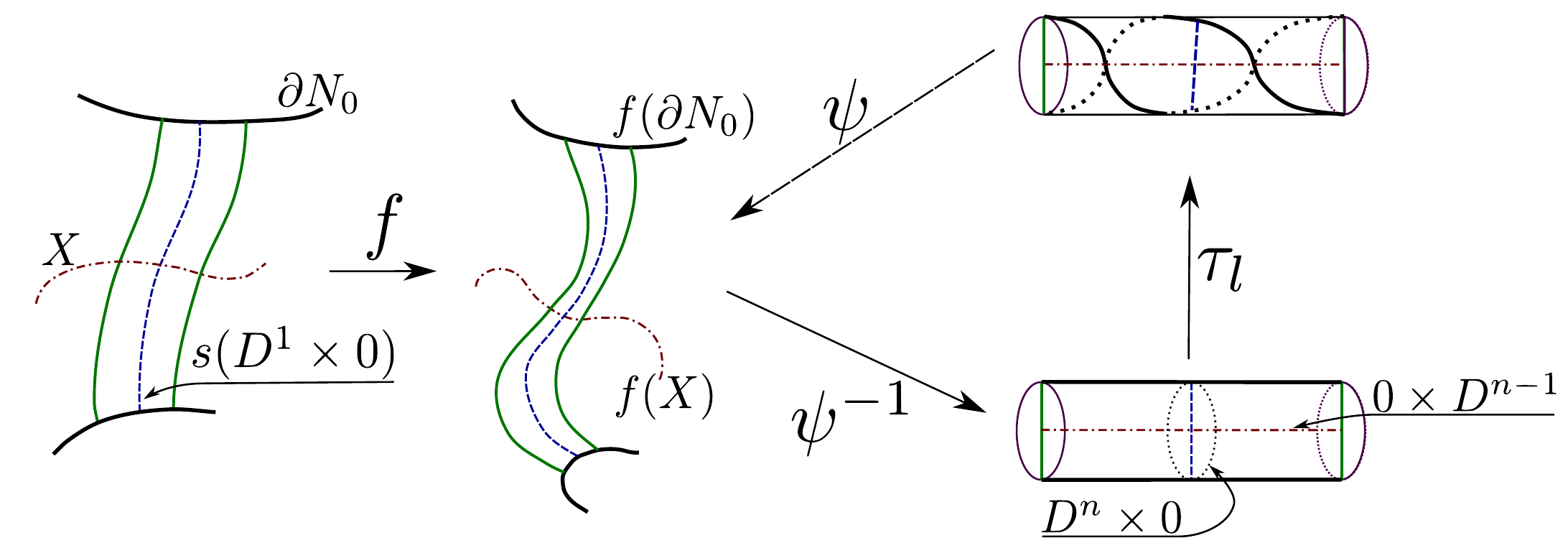}
 \label{fig::twisting}
 \caption{Twisting of a ribbon defined by $\tilde s$.}
\end{figure}

Define a map $\tau_l\colon D^1\times D^{n-1}\to D^n\times D^{n-1}$ by the formula
$$\tau_l(t, a) = (h_lq(a)\cdot t, a),$$
where $q\colon D^{n-1}\to S^{n-1}$ is the shrinking of the boundary of a disk to a point and $h_l\colon S^{n-1}\to S^{n-1}$ is a map of degree $l$, cf. Example~\ref{exm::linked_boundary}.
Hence, $\tau_l$ coincides on $D^1\times \partial D^{n-1}$ with the standard inclusion.
Then we can define an embedding
$$
g\colon N_0\to \mathbb R^{2n-1} \;\text{  by  }\; g(a):=
\begin{cases}
	f(a), & a\notin s(D^1\times D^{n-1}), \\
	\psi\tau_l \psi^{-1}f(a), & a \in s(D^1\times D^{n-1}).
\end{cases}
$$

Now let us show that
$L(g)(x, y) - L(f)(x, y) = 2l(\tilde s\cap x) \cdot(\tilde s\cap y)$
for any $x, y \in H_{n-1}(N_0)$.

Represent $x$ and $y$ by $(n-1)$-cycles $X$ and $Y$ such that
$$f(X)\cap s(D^1\times D^{n-1}) \hm= f(Y)\cap s(D^1\times D^{n-1}) = s(0\times D^{n-1}).$$
So $f(X) = g(X)$ and $f(Y) = g(Y)$.
Then
$$
\psi^{-1}(f(X)) \hm= (\tilde s\cap X)(0\times D^{n-1})
\quad\text{and }\quad
\psi^{-1}(f(Y)) = (\tilde s\cap Y)(0\times D^{n-1})
$$
are $(n-1)$-cycles in $D^n\times D^{n-1}$ modulo $D^n\times\partial D^{n-1}$.

In the following paragraph we construct a normal field $\xi_g$ on $g(N_0)$.

There exist bundle $S^{n-2}\to V_{n, 2}\stackrel{p}\to S^{n-1}$, where $p$ is the forgetful map sending $2$-frame to its first vector.
This bundle has the following exact sequence
$$
\cdots\to\pi_{n-1}(V_{n, 2})\stackrel{p_*}{\to} \pi_{n-1}(S^{n-1})\stackrel{\partial}{\to}\pi_{n-2}(S^{n-2})\to\cdots
$$
Since $n$ is even, it follows that the map $\partial$ is zero, see \cite[\S19.C]{FF16}. So the map $p_*$ is surjective.
Thus, there exist a map $h_l'\colon S^{n-1}\to V_{n, 2}$ such that $p_*h_l' \simeq h_l$.
Any element in $V_{n, 2}$ gives a map $D^2\to D^n$, so the map $h_l'$ gives a map $S^{n-1}\times D^2\to D^n$ which we denote by the same notation.
Define a map $\tau_l'\colon D^2\times D^n\to D^n\times D^{n-1}$ by the formula
$$\tau_l'(b, a) = (h_l'(q(a), b), a).$$
Since we may assume that $h_l'(q(\partial D^{n-1})) = (1_0, 1_1)\in V_{n, 2}$, it follows that $\tau_l'$ coincides on $D^2\times \partial D^{n-1}$ with the standard inclusion.
Then we can define a normal field
$$
\xi_g\colon N_0\to \mathbb R^{2n-1} \;\text{  by  }\;
\xi_g(a):=
\begin{cases}
	\xi_f(a),
	& a\notin s(D^1\times D^{n-1}), \\
	\psi\tau_l' \psi^{-1}\xi_f(a),
	& a \in s(D^1\times D^{n-1}).
\end{cases}
$$

Let $\delta f(X)$ be an $n$-chain in $\mathbb R^{2n-1}$ such that
$$
\partial\delta f(X) = f(X) = g(X)
\quad\text{ and }\quad
\psi^{-1}(\delta f(X)) \hm= [0,1_0]\times D^{n-1}
.
$$
Thus, $\delta f(X)\cap \psi(D^n\times D^{n-1}) \subset fs(D^1\times D^{n-1})$.
We have
\begin{align*}
 \lk (g(X), \xi_g(Y)) &- \lk(f(X), \xi_f(Y)) = \delta f(X) \cap \xi_g(Y) - \delta f(X) \cap \xi_f(Y)
 =\\&=
 \delta f(X) \cap (\xi_g(Y) - \xi_f(Y))
 \stackrel{(1)}{=}
 \psi^{-1}(\delta f(X))\cap\psi^{-1}(\xi_g(Y) - \xi_f(Y))
 =\\&=
 (\tilde s\cap X) ([0,1_0]\times D^{n-1}) \cap (\tilde s \cap Y) (\tau_l'(1_1\times D^{n-1}) - 1_1\times D^{n-1})
 \stackrel{(2)}{=}\\&\stackrel{(2)}{=}
 (\tilde s\cap X)\cdot(\tilde s\cap Y) ((1_0\times D^{n-1}) \cap \tau_l'(1_1\times D^{n-1}))
 \stackrel{(3)}{=}
 (\tilde s\cap X)\cdot(\tilde s\cap Y) l.
 \end{align*}

The equality $(1)$ holds because $\xi_f$ and $\xi_g$ differs only on $s(D^1\times D^{n-1})$, so $\xi_g(y) - \xi_f(y)\hm\subset \psi(D^n\times D^{n-1})$.

The equality $(2)$ holds because $[0,1_0]\times D^{n-1}$ does not intersect $1_1\times D^{n-1}$.

The equality $(3)$ holds because the set $(1_0\times D^{n-1}) \cap \tau_l'(1_1\times D^{n-1})$ consist of all points $(1_0, a)\in S^{n-1}\times S^{n-2}$ such that $1_0=h_l'(q(a), 1_1)$. The sum of their signs is $\deg h_l'(b, 1_1)$.
We have $\deg h_l'(b,1_1) = \deg h_l'(b, 1_0)\hm = \deg h_l = l$, because $h_l = h_l'(b, 1_0)$ and $h_l'(b, 1_1)$ are homotopic.

Analogously $\lk (\xi_g(X), g(Y)) - \lk(\xi_f(X), f(Y)) =  (\tilde s\cap X)\cdot(\tilde s\cap Y) l$.
This completes the proof.
\end{proof}

\section*{Acknowledgements}
The author is grateful to his advisor professor A.~Skopenkov for valuable discussions, constant support and attention to this work%
.

\end{document}